# CONFIDENCE BALLS IN GAUSSIAN REGRESSION

By Yannick Baraud

*Ecole Normale Supérieure*

Starting from the observation of an $\mathbb{R}^n$-Gaussian vector of mean $f$ and covariance matrix $\sigma^2 I_n$ ($I_n$ is the identity matrix), we propose a method for building a Euclidean confidence ball around $f$, with prescribed probability of coverage. For each $n$, we describe its nonasymptotic property and show its optimality with respect to some criteria.

**1. Introduction.** In the present paper, we consider the statistical model

$$(1) \qquad Y_i = f_i + \sigma \varepsilon_i, \qquad i = 1, \ldots, n,$$

where $f = (f_1, \ldots, f_n)'$ is an unknown vector, $\sigma$ a positive number and $\varepsilon_1, \ldots, \varepsilon_n$ a sequence of i.i.d. standard Gaussian random variables. For some $\beta \in \,]0,1[$, the aim of this paper is to build a nonasymptotic Euclidean confidence ball for $f$ with probability of coverage $1 - \beta$ from the observation of $Y = (Y_1, \ldots, Y_n)'$.

This statistical model includes, as a particular case, the functional regression model

$$(2) \qquad Y_i = F(x_i) + \sigma \varepsilon_i, \qquad i = 1, \ldots, n,$$

where $F$ is an unknown function on some interval, say $[0, 1]$, and the $x_i$'s are some distinct deterministic points in this interval. The literature on the topic usually deals with this particular model, which offers the advantage of focusing on the quantity $F$, which does not depend on $n$. This simplifies the asymptotic point of view. For this reason, we shall focus in this Introduction on the problem of building a confidence ball for $F$. In the sequel, we denote by $\|\cdot\|_n$ the seminorm defined on the set of real-valued functions $t$ on $[0, 1]$ by $\|t\|_n^2 = n^{-1} \sum_{i=1}^n t^2(x_i)$.

The problem of building a confidence ball for $F$ with respect to $\|\cdot\|_n$ easily reduces to that of building a Euclidean confidence ball for the vector $f = (F(x_1), \ldots, F(x_n))'$ by identifying the functions $t$ on $[0,1]$ with the $\mathbb{R}^n$-vectors $(t(x_1), \ldots, t(x_n))'$. Thus, when $\sigma^2$ is known, say equal to 1, the problem is solved by considering the Euclidean ball centered at $Y$ with squared radius $q_{0,n}(\beta)$, where $q_{0,n}(\beta)$ denotes the $(1-\beta)$-quantile of a $\chi^2$-distribution with $n$ degrees of freedom. However, such a confidence ball is almost useless: besides providing a very rough estimator of $F$, the radius of the confidence ball is very large. To overcome this problem, a natural idea is to start with a "good" estimator of $F$, say $\hat{F}_n$, and then to estimate $\delta_n(F) = \|F - \hat{F}_n\|_n^2$ by some suitable estimator, say $\hat{\delta}_n$. This is the key point of the procedures proposed by Li (1989), Beran (1996) and Beran and Dümbgen (1998). In the last two papers, the estimators $\hat{F}_n$ and $\hat{\delta}_n$ are such that $\sqrt{n}(\delta_n(F) - \hat{\delta}_n)$ converges to some limit distribution $Q$ as $n$ becomes large. Thus, if one denotes by $Q^{-1}(1-\beta)$ the $(1-\beta)$-quantile of $Q$, the ball centered at $\hat{F}_n$ of squared radius $\hat{\delta}_n + Q^{-1}(1-\beta)/\sqrt{n}$ provides a confidence region with asymptotic probability of coverage $1-\beta$. The limit distributions $Q$ obtained in Beran (1996) and Beran and Dümbgen (1998) are both Gaussian of mean 0. However, their variances depend on $F$ and $\sigma$ and, consequently, $Q^{-1}(1-\beta)$ must be estimated in turn from the data. The disadvantage of the procedures proposed in Beran (1996) and Beran and Dümbgen (1998) mainly lies in their asymptotic character. It is indeed difficult to judge whether the asymptotic regime is achieved or not as it depends on the features of the unknown function $F$.

In contrast, the asymptotic confidence balls proposed by Li (1989) are called *honest* in the sense that the probability of coverage is uniform with respect to all possible functions $F$. However, in Li (1989) the variance of the errors is assumed to be known and the radius of the confidence ball involves an inexplicit constant. His procedure is based on a Stein estimator of $F$, $\hat{F}_n$, and a Stein estimator of $\|F - \hat{F}_n\|_n^2$. A comparison between Li's confidence balls and ours will be given in Section 2.3.

Another direction was investigated by Cox (1993). He considered Bayesian inference for a class of regression models. The regression functions $F$ were drawn under a Gaussian prior distribution among the solutions of a high-order stochastic differential equation. He analyzed the $\mathbb{L}^2([0,1], dx)$-distance between $F$ and its estimator $\hat{F}$ (the posterior expectation of $F$) and deduced a confidence ball for $F$. He proved that if $n$ is fixed (large enough) the frequentist probability of coverage of the confidence ball is close to 1 for all $F$ within a set of probability close to 1. However, this probability of coverage is infinitely often less than any positive $\varepsilon$ as $n$ tends to infinity for almost all $F$. Unfortunately, this negative result on Cox's confidence ball makes it unattractive for non-Bayesians.

The ideas underlying our approach are due to Lepski and have been exposed by their initiator in a series of lectures at the Institute Henri Poincaré in Paris. We shall now give a brief account of these ideas and recommend that the reader have a look at Lepski (1999) for more details. Lepski noted that if $F$ is known to belong to a suitable class $\Sigma$ of smooth functions, then the minimax approach allows one to obtain both an estimator of $F$ and a control on the accuracy of the estimation. However, unless one has a strong guess on the particular features of $F$, $\Sigma$ is usually too large to obtain an accurate estimation. The idea of Lepski is to test one or several additional structures on $F$ in order to improve the accuracy of estimation. Unlike an adaptive approach, an attractive feature of Lepski's approach lies in that the accuracy is available to the statistician and, consequently, that a nonparametric confidence ball for $F$ can be derived. This is explained in the papers by Lepski (1999) and by Hoffmann and Lepski (2002). However, the procedure described there for the purpose of building $\mathbb{L}^2$-confidence balls suffers from the following weaknesses. First, the point of view is purely asymptotic. The procedure does not lead to confidence balls with prescribed probability of coverage for fixed values of $n$. Furthermore, a careful look at the proofs shows that, for a fixed $n$, the squared radius of the confidence ball is equal to a constant plus some term which is essentially proportional to the number of hypotheses to test. Consequently, the number of these cannot be large if one wants to keep the confidence ball of a reasonable size. In addition, the squared radius of the confidence ball is proportional to $1/\beta$ and is thus very large for small values of $\beta$. Finally, the applications developed in Lepski (1999) and Hoffmann and Lepski (2002) mainly address the Gaussian white noise model and an adaptation of the procedure to the regression case would require an estimation of the unknown $\sigma$.

The results of the present paper are nonasymptotic and the procedures which are described here aim at obtaining confidence balls which are as sharp as possible. In particular, the dependency with respect to $\beta$ and the number of hypotheses to test is only logarithmic. This allows us to handle the variable selection problem described in Section 2.4.

We consider the case where $\sigma$ is known to belong to some interval $I = [(1-\eta)\tau^2, \tau^2]$ with $\eta \geq 0$. The situation $\eta = 0$ corresponds to the theoretical situation where one exactly knows the variance. In contrast, the situation $\eta > 0$ corresponds to the practical one when the variance is known to belong to some interval which is either derived by the experimental context or by statistical estimation (from an independent sample). In all cases, the optimality (in a suitable sense) of our confidence balls is established. The proof relies on nonasymptotic lower bounds for the minimax estimation and separation rates over linear spaces. We show that if a confidence ball ensures the probability of coverage $1 - \beta$ uniformly over all $f \in \mathbb{R}^n$ and $\sigma^2 \in I$, then its radius (normalized by $\sqrt{n}$) must be greater than $C \max\{\sqrt{\eta}, n^{-1/4}\}$, where $C$ is a constant free from $n$ and $\eta$. When $\eta = 0$, this result allows one to recover that established by Li (1989), namely that asymptotically the radius of such a confidence ball cannot converge toward 0 faster than $n^{-1/4}$. When $\eta > 0$, this result shows that practically the problem of establishing useful confidence balls is impossible unless $\eta$ is small compared to $n$.

The paper is organized as follows. In Section 2 we consider the case of a known $\sigma(\eta = 0)$ and describe a procedure free from any prior assumption on $f$. This procedure is implemented on numerical examples in Section 4. In Section 3, we consider the case $\eta > 0$ and provide some lower bounds on the radius of an honest confidence ball. We show in this section that these lower bounds are sharp by providing a construction of confidence balls which achieves these bounds. The proofs are postponed to Section 5.

NOTATION.  Throughout this paper we use the following notation. We denote by $\|\cdot\|$ the Euclidean distance in $\mathbb{R}^n$. For a triplet $(z, d, u) \in \mathbb{R}_+ \times \mathbb{N} \setminus \{0\} \times {]0, 1[}$, we denote by $\chi^2_{z,d}(\cdot)$ the distribution function of a (non)central $\chi^2$ with noncentrality parameter $z$ and $d$ degrees of freedom and by $q_{z,d}(u)$ its $(1-u)$-quantile for $u \in {]0, 1[}$. In particular, if $X$ is distributed as $\chi^2_{z,d}(\cdot)$, then

$$\mathbb{E}[X] = z + d, \quad \text{and} \quad \mathbb{P}(X \geq q_{z,d}(u)) = u \qquad \forall u \in {]0, 1[}.$$

We will use the convention $q_{z,0}(u) = 0$ for all $u \in {]0,1[}$ and $z \geq 0$. For each linear subspace $S$ of $\mathbb{R}^n$, we denote by $\Pi_S$ the orthogonal projector onto $S$ and by $\mathcal{B}(x, r)$ the Euclidean ball centered at $x \in \mathbb{R}^n$ of radius $r > 0$. Finally, $C, C', \ldots$ denote constants that may vary from line to line.

**2. Confidence balls when the variance is known.**  The aim of this section is twofold: first, explain the basic ideas of our approach and second, in the ideal case where the variance $\sigma^2$ is known, build a confidence ball for $f$ with controlled probability of coverage.

2.1. *The basic ideas.*  An ideal procedure to build a confidence ball would probably be to start with a nice estimator of $f$, say $\hat{f}$, and then get a uniform control of $\|f - \hat{f}\|$ over all possible $f$. This strategy is unfortunately impossible in general. For illustration, let us consider $\hat{f} = \Pi_S Y$, the projection estimator of $f$ onto a linear subspace $S$ of $\mathbb{R}^n$ of dimension $\mathcal{D} < n$. By setting $z$ equal to the squared Euclidean distance between $f$ and $S$ and using Pythagoras' theorem, we derive that

$$\|f - \hat{f}\|^2 = z + \|\Pi_S \varepsilon\|^2 \sigma^2$$

and, hence, a control of $\|f - \hat{f}\|^2$ necessarily requires that an upper bound on $z$ be known. This is of course seldom the case in practice. The idea of our procedure is to get such a piece of information by means of a test. More precisely, let us fix some $\alpha \in ]0, 1-\beta[$ and consider the $\chi^2$-test of level $\alpha$ of hypothesis "$f \in S$" against "$f \in \mathbb{R}^n \setminus S$" which consists in rejecting the null when the test statistic $T = \|Y - \Pi_S Y\|^2$ is greater than $q_{0,n-\mathcal{D}}(\alpha)\sigma^2$. If the test accepts the null, then intuitively this means that $f$ is close to $S$ and, therefore, that $z$ is small. The following lemma shows that $\|f - \hat{f}\|$ cannot be large on the event that the hypothesis "$f \in S$" is accepted.

LEMMA 2.1. *Let $\alpha \in ]0, 1-\beta[$. Let us define*

$$\phi(Y) = \mathbf{1}\{\|Y - \Pi_S Y\|^2 > q_{0,n-\mathcal{D}}(\alpha)\sigma^2\} \tag{3}$$

*and*

$$\mathcal{Z} = \{z \in \mathbb{R}_+, \chi^2_{z,n-\mathcal{D}}(q_{0,n-\mathcal{D}}(\alpha)) > \beta\}.$$

*If $\mathcal{D} \neq 0$, we set*

$$\rho^2 = \sup_{z \in \mathcal{Z}} \left[ z + q_{0,\mathcal{D}}\left( \frac{\beta}{\chi^2_{z,n-\mathcal{D}}(q_{0,n-\mathcal{D}}(\alpha))} \right) \right] \sigma^2; \tag{4}$$

*if $\mathcal{D} = 0$, we set*

$$\rho^2 = \inf\{z \geq 0, \chi^2_{z,n}(q_{0,n}(\alpha)) \leq \beta\}\sigma^2. \tag{5}$$

*Then, for all $f \in \mathbb{R}^n$,*

$$\mathbb{P}_{f,\sigma}[\phi(Y) = 0, \|f - \hat{f}\| \geq \rho] \leq \beta. \tag{6}$$

Let us assume that $\sigma = 1$ and make a few comments on the set $\mathcal{Z}$ and the quantity $\rho$. The inequality $\alpha < 1 - \beta$ implies that 0 belongs to $\mathcal{Z}$ and, hence, the set $\mathcal{Z}$ is always nonvoid. Moreover, since the map $\psi : z \mapsto \chi^2_{z,n-\mathcal{D}}(q_{0,n-\mathcal{D}}(\alpha))$ is decreasing, continuous and tends to 0 as $z$ becomes large, it appears that $\mathcal{Z}$ is an interval of the form $[0, \bar{z}[$, where $\bar{z}$ satisfies $\psi(\bar{z}) = \beta$. When $\mathcal{D} = 0$ we deduce that $\rho^2 = \bar{z}$ and, consequently, that $\rho$ is finite. Since $q_{0,\mathcal{D}}(u)$ tends to 0 as $u$ approaches 1 from below, we see that $\rho^2$ is also finite when $\mathcal{D} \neq 0$. The supremum in (4) is usually achieved at some point $z^* \in \mathcal{Z}$. If the squared Euclidean distance between $f$ and $S$ equals $z^*$, then equality holds in (6). The quantity $z^*$ is a critical value for the (squared) distance $z$ between $f$ and $S$: if $z$ is large compared to $z^*$, then the test $\phi$ rejects the null with probability close to 1 and thus the left-hand side of (6) is small. This is also the case if, on the other hand, $z$ is small compared to $z^*$ because then $\hat{f}$ is a "good" estimator of $f$ and the event $\|f - \hat{f}\| > \rho$ seldom occurs.

The convention

$$q_{0,\mathcal{D}}(1) = -\infty \tag{7}$$

allows one to define the quantity $\rho$ equivalently as

$$\tag{8} \rho^2 = \sup_{z \geq 0} \left[ z + q_{0,\mathcal{D}} \left( \frac{\beta}{\chi^2_{z,n-\mathcal{D}}(q_{0,n-\mathcal{D}}(\alpha))} \wedge 1 \right) \right] \sigma^2.$$

In the sequel, we shall use this convention to simplify our notation.

Our procedure for building a confidence ball around $f$ is based on Lemma 2.1. As a control of $\|f - \hat{f}\|$ is possible when the hypothesis "$f \in S$" is accepted, we increase our chance to accept such hypotheses by considering a family of $S$'s rather than a single one. Moreover, in order to ensure that, for at least one $S$ the hypothesis "$f \in S$" is accepted, we add the linear space $S = \mathbb{R}^n$ to the family, the hypothesis "$f \in \mathbb{R}^n$" being obviously true.

2.2. *Construction of the confidence ball.* Let $\{S_m, m \in \mathcal{M}_n\}$ be a finite family of linear subspaces of $\mathbb{R}^n$. For each $m$, we set $\mathcal{D}_m = \dim(S_m)$, $N_m = n - \mathcal{D}_m$ and associate with $S_m$ some number $\beta_m$ in $]0,1[$. We assume that the following assumption is fulfilled.

ASSUMPTION 2.1. *The subscript $n$ belongs to $\mathcal{M}_n$ and $S_n = \mathbb{R}^n$. We have $\sum_{m \in \mathcal{M}_n} \beta_m \leq \beta$.*

For each $m \in \mathcal{M}_n$, we define $\rho_m$ as follows. If $m = n$, then

$$\rho_n^2 = q_{0,n}(\beta_n)\sigma^2.$$

If $m \in \mathcal{M}_n \setminus \{n\}$ and $\mathcal{D}_m \neq 0$, then $\rho_m$ is defined by (8) with $\mathcal{D}_m$ in place of $\mathcal{D}$ and $\beta_m$ in place of $\beta$. If $m \in \mathcal{M}_n \setminus \{n\}$ and $\mathcal{D}_m = 0$, then $\rho_m$ is defined by (5) with $\beta_m$ in place of $\beta$.

For each $m \in \mathcal{M}_n \setminus \{n\}$, we define $\hat{f}_m = \Pi_{S_m} Y$ and $\phi_m$ is the test defined by (3) with $S = S_m$. If $m = n$, then $\hat{f}_n = Y$ and $\phi_n(y) = 0$ for all $y \in \mathbb{R}^n$.

We define

$$\mathcal{A} = \{m \in \mathcal{M}_n, \phi_m(Y) = 0\}$$

and

$$\tag{9} \hat{m} = \arg\min_{m \in \mathcal{A}} \rho_m, \qquad \hat{\rho} = \rho_{\hat{m}}, \qquad \hat{f} = \hat{f}_{\hat{m}}.$$

We have the following result.

THEOREM 2.1. *Let $(\hat{f}, \hat{\rho})$ be the pair of random variables defined by (9). The region $\mathcal{B}(\hat{f}, \hat{\rho})$ is a confidence ball with probability of coverage $1 - \beta$, that is,*

$$\tag{10} \mathbb{P}_{f,\sigma}[f \in \mathcal{B}(\hat{f}, \hat{\rho})] \geq 1 - \beta \qquad \forall f \in \mathbb{R}^n.$$

*Moreover, for each $m \in \mathcal{M}_n$ and $f \in \mathbb{R}^n$, if for some $\gamma \in \,]0,1[$ we have*

(11) $\quad\quad \mathbb{P}_{f,\sigma}[\phi_m(Y) = 0] \geq 1 - \gamma \quad\quad \text{then } \mathbb{P}_{f,\sigma}[\hat{\rho} \leq \rho_m] \geq 1 - \gamma.$

*In particular, for all $m \in \mathcal{M}_n$,*

(12) $\quad\quad \inf_{f \in S_m} \mathbb{P}_{f,\sigma}[\hat{\rho} \leq \rho_m] \geq 1 - \alpha.$

Let us make a few comments:

1. Inequalities (11) and (12) are clear from the definition of $\hat{\rho}$ since with probability not less than $1 - \gamma$ (resp. $1 - \alpha$) we have $m \in \mathcal{A}$. Inequality (12) provides an upper bound (in probability) for the random variable $\hat{\rho}$ under the law $\mathbb{P}_{f,\sigma}$ as soon as $f \in S_m$. Inequality (11) says that this upper bound remains valid not only when $f$ belongs to $S_m$ but also when $f$ is close to $S_m$, as then the test $\phi_m$ still accepts the hypothesis "$f \in S_m$" with large probability.
2. Note that $\mathcal{A}$ is nonvoid since $n$ belongs to $\mathcal{A}$. The case where $\hat{\rho} = \rho_n$ corresponds to the one where none of the hypotheses "$f \in S_m$" (with $m \in \mathcal{M}_n \setminus \{n\}$) is accepted. In this case, the resulting confidence ball is crude, namely centered at $Y$ of radius $\rho_n$. Note that when $\beta_n$ is chosen to be of order $\beta$, say $\beta/2$, the radius $\rho_n^2$ is of the same order as $\bar{\rho}^2 = q_{0,n}(\beta)\sigma^2$, which means that the procedure does not lose too much compared with the trivial confidence ball $\mathcal{B}(Y, \bar{\rho})$.
3. In the proofs we show something stronger than Theorem 2.1. Namely, we prove that, with probability not less than $1 - \beta$, $f$ belongs to the intersection of the Euclidean balls $\mathcal{B}(\hat{f}_m, \rho_m)$ for $m \in \mathcal{A}$. However, the resulting confidence region is no longer a ball in general.

The expressions of the quantities $\rho_m$ do not allow a direct appreciation of their orders of magnitude. An upper bound for $\rho_m$ is given in the following proposition. We restrict ourselves to the case where the dimension of $S_m$ is not larger than $n/2$. Indeed, considering linear spaces with dimension larger than $n/2$ leads to large radii and thus does not offer a real gain compared to $\mathbb{R}^n$. The proof of the following proposition contains explicit constants.

PROPOSITION 2.1. *Assume that, for all $m \in \mathcal{M}_n \setminus \{n\}$, $\mathcal{D}_m \leq n/2$. Then there exists some constant $C$ depending on $\alpha$ only such that, for all $m \in \mathcal{M}_n$,*

$$\rho_m^2 \leq C \max\{\mathcal{D}_m, \sqrt{n \log(1/\beta_m)}, \log(1/\beta_m)\}\sigma^2.$$

If $\mathcal{M}_n$ reduces to $\{n\}$, then $\hat{\rho} = \rho_n$ and the radius of the ball is of order $n\sigma^2$ by taking $\beta_n = \beta$. By considering several linear spaces $S_m$ we have the opportunity to capture some specific features of $f$ and consequently to reduce the order of magnitude of $\hat{\rho}$. The number of tests $|\mathcal{M}_n|$ to perform is taken into account via the quantity $\beta_m$. If one chooses $\beta_m = \beta/|\mathcal{M}_n|$ for all $m \in \mathcal{M}_n$, one gets that the radius of the confidence ball depends logarithmically on $|\mathcal{M}_n|$. However, a choice of $\beta_m$ depending on $m$ via the dimension of the linear space $S_m$, for example, is recommended. We shall see an example in Section 2.4.

2.3. *Comparison with the procedure proposed by Li.* In this section, we make a comparison between our procedure and that proposed by Li. To simplify the discussion we assume that $\sigma^2 = 1$. Li's procedure relies on a Stein estimator of $f$, say $\tilde{f}^*$, and a Stein estimator of $\|f - \tilde{f}^*\|^2$. The estimator $\tilde{f}^*$ is obtained by modifying a linear estimator of $f$, say $\hat{f}$. By taking $\hat{f} = \Pi_S Y$, where $S$ is a linear subspace of $\mathbb{R}^n$ of dimension $\mathcal{D} < n$, the confidence ball Li proposes is centered at

$$\tilde{f}^* = \hat{f} + \left(1 - \frac{n - \mathcal{D}}{\|Y - \Pi_S Y\|^2}\right)(Y - \Pi_S Y)$$

and its squared radius is given by

$$r^2 = c\sqrt{n} + n\left(1 - \frac{(n - \mathcal{D})^2}{n\|Y - \Pi_S Y\|^2}\right),$$

where $c$ is an unspecified constant depending on $\beta$ and $\sigma^2$ only. He proved this confidence ball has probability of coverage $1 - \beta$ for all $f \in \mathbb{R}^n$ simultaneously provided that $n$ is large enough. To compare this confidence ball to ours, let us make the a posteriori assumption that $f$ belongs to $S$. On the one hand, by using our procedure with $\mathcal{M}_n = \{m, n\}$, $S_m = S$, $\beta_m = \beta/2 = \beta_n$, we derive from Theorem 2.1 that, with probability close to 1, $\hat{\rho}^2 = \rho_m^2$, which is of order $\max\{\sqrt{n}, \mathcal{D}\}$. On the other hand, replacing $\|Y - \Pi_S Y\|^2$ by its expectation $n - \mathcal{D}$ shows that the squared radius of Li's confidence ball is of order

$$r^2 \approx c\sqrt{n} + n\left(1 - \frac{n - \mathcal{D}}{n}\right) = c\sqrt{n} + \mathcal{D}$$

and is therefore of the same order as ours.

However, for those $f$ which do not belong to $S$ the radius of Li's confidence ball can become large. The advantage of our approach lies in that it is possible to deal with a larger family of spaces than just $\{S, \mathbb{R}^n\}$. By doing so, we can keep the radius of the confidence ball to a reasonable size for those vectors $f$ which are close to at least one of the linear spaces of the family and not only $S$.

2.4. *Application to variable selection.* In this section, we illustrate the procedure in the variable selection problem. Assume that $f$ is of the form $XU$, where $X$ is a known $p \times n$ full-rank matrix with $p \in \{1, \ldots, n\}$ and $U$ some unknown vector in $\mathbb{R}^p$. The problem of variable selection is to determine from the data the nonzero coordinates of $U$, that is,

$$m^* = \{j \in \{1,\ldots,p\}, U_j \neq 0\}.$$

In this section we give a way to select those coefficients and provide simultaneously a confidence ball for $f$. We apply the procedure as follows:

Let $\mathbf{x}_1,\ldots,\mathbf{x}_p$ be the column vectors of the matrix $X$ and let $\mathcal{P}_n$ be the class of nonempty subsets $m$ of $\{1,\ldots,p\}$ with cardinality $|m|$ not larger than $n/2$. For all $m \in \mathcal{P}_n$, we define $S_m$ as the linear span of the $\mathbf{x}_j$'s for $j \in m$ and set

$$\beta_m = \beta \left[ n \binom{n}{\mathcal{D}} \right]^{-1} \qquad \text{with } \mathcal{D} = |m|.$$

We define $\mathcal{M}_n = \mathcal{P}_n \cup \{n\}$ and set $\beta_n = \beta/2$. Note that Assumption 2.1 is fulfilled since

$$\sum_{m \in \mathcal{M}_n} \beta_m = \frac{\beta}{2} + \sum_{m \in \mathcal{P}_n} \beta_m = \frac{\beta}{2} + \sum_{1 \leq \mathcal{D} \leq n/2} \sum_{m \in \mathcal{P}_n, |m|=\mathcal{D}} \beta_m \leq \beta.$$

By applying the procedure described in Section 2.2 we select a set of indices $\hat{m}$ for which the Euclidean distance between the least-squares estimator $\hat{f}_{\hat{m}}$ and $f$ is not greater than $\rho_{\hat{m}}$ with probability greater than $1 - \beta$. Since $f$ belongs to the linear space $S_{m^*}$, with probability greater than $1 - \alpha$ the set $m^*$ belongs to $\mathcal{A}$ and consequently $\rho_{\hat{m}}$ is not greater than $\rho_{m^*}$. Therefore, either $\hat{m} = m^*$ and then the procedure selects the target subset $m^*$, or $\hat{m} \neq m^*$ and then the resulting confidence ball is at least as accurate as if the target subset $m^*$ were selected. In addition, thanks to the inequality

$$\binom{n}{D} \leq \exp(\mathcal{D} \log(en/\mathcal{D}))$$

and Proposition 2.1, with probability greater than $1 - \alpha$, the following upper bound holds: there exists some constant $C$ depending on $\alpha$ and $\beta$ only such that

$$\hat{\rho}^2 \leq C \max\{\sqrt{n|m^*|\log(en/|m^*|)}, |m^*|\log(en/|m^*|)\}\sigma^2.$$

Let us denote this upper bound by $B$. Another possible choice of the $\beta_m$'s is $\beta_m = \beta_n = \beta/|\mathcal{M}_n|$ for all $m \in \mathcal{M}_n$. For this second strategy, $\hat{\rho}^2$ is of order $B' = \max\{\sqrt{np}, p\}\sigma^2$ as $|\mathcal{M}_n|$ is of order $2^p$. In the least favorable situation where almost all the coefficients $U_j$'s are nonzero, $|m^*|$, $p$ and $n$ are of the same order and, thus so are $B$ and $B'$. In this case, both strategies lead to confidence balls which are approximately of the same size. Yet, in the more favorable situation where $p$ is still of order $n$ but $|m^*|$ is small compared to $p$, the strategy with nonconstant $\beta_m$'s leads to a sharper confidence ball. This illustrates the advantage of taking $\beta_m$ as a function of $m$.

**3. Confidence balls under some information on the variance.** In this section, we no longer assume that $\sigma$ is known but rather that it belongs to some known interval $I = [\sqrt{1-\eta}\tau, \tau]$, where $(\tau^2, \eta) \in \mathbb{R}_+ \times [0,1[$. As we shall see, the uncertainty on the value of $\sigma$ has a terrible effect on the orders of magnitude of radii of confidence balls.

3.1. *How sharp can the confidence ball be?* We have the following result.

THEOREM 3.1. *Let $\alpha$ and $\beta$ be numbers in $]0,1[$ satisfying $2\beta + \alpha < 1 - \exp(-1/36)$. Let $(\tilde{f}, \tilde{r})$ be a pair of random variables depending on $Y$ only with values in $\mathbb{R}^n \times \mathbb{R}_+$ satisfying, for all $f \in \mathbb{R}^n$ and $\sigma \in I$,*

(13) $$\mathbb{P}_{f,\sigma}[f \in \mathcal{B}(\tilde{f}, \tilde{r})] \geq 1 - \beta.$$

*For each $m \in \mathcal{M}_n$, let $r_m$ be some positive quantity satisfying for all $\sigma \in I$*

(14) $$\inf_{f \in S_m} \mathbb{P}_{f,\sigma}[\tilde{r} \leq r_m] \geq 1 - \alpha.$$

*Then there exists some constant $C$ depending on $\alpha$ and $\beta$ only such that, for all $m \in \mathcal{M}_n$,*

(15) $$r_m^2 \geq C \max\{\eta N_m, \mathcal{D}_m, \sqrt{N_m}\}\tau^2.$$

*For each $f \in \mathbb{R}^n$ let $r(\alpha, f)$ be such that, for all $\sigma \in I$,*

$$\mathbb{P}_{f,\sigma}\{\tilde{r} \leq r(\alpha, f)\} \geq 1 - \alpha.$$

*Then we have*

(16) $$r^2(\alpha, f) \geq C \max\{\eta n, \sqrt{n}\}\tau^2.$$

To keep our formula as legible as possible, the above theorem involves an inexplicit constant $C$. However, lower bounds including explicit numerical constants are available from the proof in Section 5.3.

Let us make few comments.

1. From an asymptotic point of view, (16) allows one to recover the result established by Li, namely that the radius of an honest confidence ball (normalized by $\sqrt{n}$) cannot converge toward 0 faster than $n^{-1/4}$. We also get that the thus normalized radius converges towards 0 only if $\eta = \eta(n)$ does and then the rate cannot be better than $\max\{\sqrt{\eta(n)}, n^{-1/4}\}$.
2. When $\eta = 0$ and $\mathcal{D}_m \leq n/2$ we derive from (15) that

   $$r_m^2 \geq C \max\{\mathcal{D}_m, \sqrt{n}\}\sigma^2,$$

   for some constant $C$ depending on $\alpha$ and $\beta$ only. This lower bound is of the same order as the upper bound on $\rho_m^2$ established in Proposition 2.1 provided that $\beta_m$ is free from $n$. This is the case if $\beta_m = \beta/|\mathcal{M}_n|$ and if the cardinality of the collection, $|\mathcal{M}_n|$, does not depend on $n$. The procedure is then optimal in the sense given by Lepski (1999).

A natural idea to establish a confidence ball around $f$ when the true variance is unknown is to use the construction of the previous section and to replace the variance $\sigma$ by the upper bound $\tau$, this latter quantity being connected "intuitively" to the least favorable situation where the level of the noise is maximal. Unfortunately, Theorem 3.1 says that such a construction cannot lead to a confidence ball as changing $\sigma$ into $\tau$ would only affect the order of magnitude of the radius by a factor $\tau/\sigma$, which would be contradictory with (16). In the next section, we show how to modify our previous construction (with a known $\sigma$) in view of obtaining a confidence ball whatever the values of $f$ and $\sigma \in I$.

3.2. *Construction of a confidence ball.* In this section we build a confidence ball under the information that $\sigma$ belongs to $I$.

The following result holds.

THEOREM 3.2. *Let $\sigma \in I$ and assume that Assumption 2.1 is fulfilled. Consider the construction of $(\hat{f}, \hat{\rho})$ described in Section 2.2 with the following definitions for the $\rho_m$'s and $\mathcal{A}$: if $m = n$, then*

$$\rho_n^2 = q_{0,n}(\beta_n)\tau^2;$$

*if $m \in \mathcal{M}_n \setminus \{n\}$ and $\mathcal{D}_m \neq 0$,*

$$\rho_m^2 = \sup_{z \geq 0, \sigma \in I} \left[ z\sigma^2 + q_{0,\mathcal{D}_m}\left( \frac{\beta_m}{\chi^2_{z,N_m}(q_{0,N_m}(\alpha)\tau^2/(\sigma^2))} \wedge 1 \right)\sigma^2 \right];$$

*if $m \in \mathcal{M}_n \setminus \{n\}$ and $\mathcal{D}_m = 0$,*

$$\rho_m^2 = \inf\left\{ x \geq 0, \sup_{\sigma \in I} \chi^2_{x/\sigma^2, n}(q_{0,n}(\alpha)\tau^2/\sigma^2) \leq \beta_m \right\}$$

*and*

$$\mathcal{A} = \{ m \in \mathcal{M}_n, \|Y - \hat{f}_m\|^2 \leq q_{0,N_m}(\alpha)\tau^2 \}.$$

*The region $\mathcal{B}(\hat{f}, \hat{\rho})$ is a confidence ball with probability of coverage $1 - \beta$; that is, (10) is satisfied. Moreover, for each $m \in \mathcal{M}_n$,*

(17) $$\inf_{f \in S_m} \mathbb{P}_{f,\sigma}[\hat{\rho} \leq \rho_m] \geq 1 - \alpha.$$

An upper bound for $\rho_m$ is given by the following proposition.

PROPOSITION 3.1. *Assume that, for all $m \in \mathcal{M}_n \setminus \{n\}$, $\mathcal{D}_m \leq n/2$. There exists some constant $C$ depending on $\alpha$ only such that, for all $m \in \mathcal{M}_n$,*

$$\rho_m^2 \leq C \max\{\eta n, \mathcal{D}_m, \sqrt{n \log(1/\beta_m)}, \log(1/\beta_m)\}\tau^2.$$

From an asymptotic point of view, we derive from Theorem 3.1 the optimality of the procedure whenever the cardinality of the collection $|\mathcal{M}_n|$ does not depend on $n$ by taking $\beta_m = \beta/|\mathcal{M}_n|$ for all $m \in \mathcal{M}_n$. For more general collections, the procedure is also optimal for those $m \in \mathcal{M}_n$ for which $\beta_m$ does not decrease with $n$.

**4. Illustrative numerical examples.** In this section we apply our procedure in three examples. In the sequel, the number of observations is $n = 1000$. We choose $\beta = 10\%$ and $\alpha = 20\%$. The $\varepsilon_i$'s are standard i.i.d. Gaussian random variables and we assume that the variance is known, that is, $\sigma^2 = 1$. We set $x_i = i/n$ for $i = 1, \ldots, n$ and define the vector $f$ as $(F(x_1), \ldots, F(x_n))'$, where $F$ is one of the following functions on $[0, 1]$:

$$F_1(x) = \cos(2\pi x),$$
$$F_2(x) = \cos(2\pi x) + 0.3 \sin(20\pi x),$$
$$F_3(x) = \begin{cases} 1.5, & \text{if } 0 < x < 0.3, \\ 0.5, & \text{if } 0.3 < x < 0.6, \\ 2, & \text{if } 0.6 < x < 0.8, \\ 0, & \text{else.} \end{cases}$$

For each function $F \in \{F_1, F_2, F_3\}$, Figure 1 shows $F$ with one set of simulated data.

For each $m \geq 1$, we define $\mathcal{F}_m$ as the linear span generated by the constant function on $[0, 1]$, $\phi_0 \equiv 1$, together with the sine and cosine functions $\cos(2\pi j x), \sin(2\pi j x)$ for $j = 1, \ldots, m$. For each $m \geq 1$, we define $S_m$ as the linear space

$$S_m = \{(F(x_1), \ldots, F(x_n))', F \in \mathcal{F}_m\}.$$

We take

$$\mathcal{M}_n = \{2^k, k = 1, \ldots, K_n\} \cup \{n\},$$

with $K_n = 8$. The number $K_n$ is chosen such that $\dim(S_{2^{K_n}}) < n$. We choose $\beta_n = \beta 2^{-K_n}$ and for each $k = 1, \ldots, K_n$, $\beta_{2^k} = \beta 2^{-k}$.

We made 100 simulations. For each simulation and each function $F \in \{F_1, F_2, F_3\}$ we consider $m(F)$, the smallest integer $m \in \mathcal{M}_n$ such that the hypothesis "$f \in S_m$" is accepted. In Table 1 we have displayed for each $F$ and $m \in \mathcal{M}_n$ the number of simulations for which $m(F) = m$.

Let us now comment on Table 1. Note that the radii $\rho_m$'s are increasing with $\mathcal{D}_m$. This comes from our choices of $\beta_m$'s, which are more favorable to linear spaces with small dimensions. Thus, the smaller is the dimension $S_m$, the sharper is the radius of the confidence ball when the hypothesis "$f \in S_m$" is accepted.

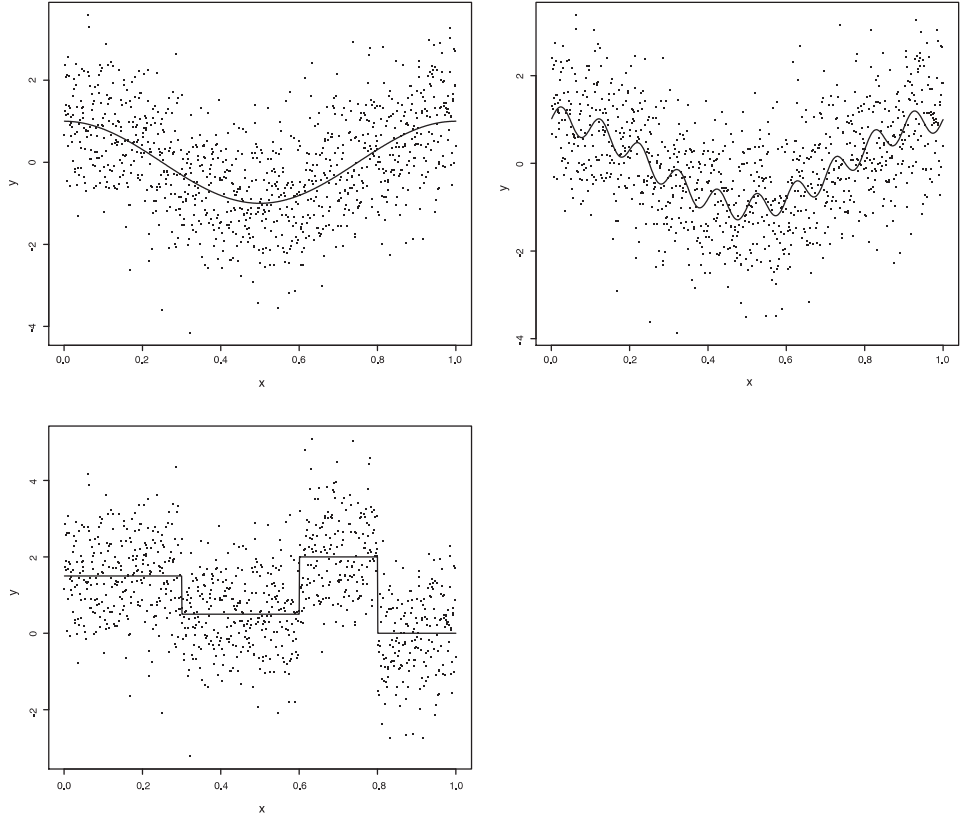

Fig. 1.

Table 1

| Indices $m$ | Dimensions $\mathcal{D}_m$ | Squared radii $\rho_m^2/n$ | "$f \in S_m$" | | |
|---:|---:|---:|---:|---:|---:|
| | | | $F_1$ | $F_2$ | $F_3$ |
| 2 | 5 | 0.118 | 82 | 47 | 0 |
| 4 | 9 | 0.136 | 1 | 0 | 8 |
| 8 | 17 | 0.155 | 0 | 1 | 20 |
| 16 | 33 | 0.181 | 1 | 33 | 28 |
| 32 | 65 | 0.222 | 1 | 3 | 17 |
| 64 | 129 | 0.293 | 4 | 5 | 6 |
| 128 | 257 | 0.425 | 1 | 1 | 7 |
| 256 | 513 | 0.681 | 4 | 4 | 5 |
| 1000 | 1000 | 1.157 | 6 | 6 | 9 |

Function $F_1$ belongs to $\mathcal{F}_2$. As expected, the hypothesis "$f \in S_2$" is accepted for around 80 simulations, $\alpha = 20\%$. This choice of $\alpha$ is arbitrary. By taking $\alpha$ smaller, the hypothesis will be accepted more often but on the other hand the radius of the confidence ball will be larger. For example, the value of $\rho_2^2/n$, respectively, equals 0.149 and 0.160 for $\alpha = 15\%$ and $\alpha = 10\%$.

Function $F_2$ is a perturbation of $F_1$. The test "$f \in S_2$" is accepted for 47 simulations even though $F_2$ does not belong to $\mathcal{F}_2$ but $\mathcal{F}_{16}$. However, for these 47 simulations the procedure has taken advantage of the closeness of $F_2$ to $\mathcal{F}_2$ to provide a sharper confidence ball than the one we would obtain if $m(F_2)$ were equal to 16. We emphasize that the procedure provides a confidence ball with probability of coverage 90% even though the "right" model for $F_2$ (namely $\mathcal{F}_{16}$) is accepted for only 33 simulations. This comes from the fact that the radius of the confidence ball takes into account a possible bias between the true and the linear space accepted by the test. Finally note that, as expected from Theorem 2.1, the radius of the confidence ball exceeds $\rho_{16}^2/n$ for 19 simulations since $F_2$ belongs to $\mathcal{F}_{16}$.

Function $F_3$ was considered in Beran and Dümbgen (1998) in one simulated example. In their simulation, the squared radius (with respect to $\|\cdot\|/\sqrt{n}$) of the confidence ball was obtained by bootstrap and was equal to 0.144. We obtain a radius of the same order for $28 = 8 + 20$ simulations.

**5. Proofs.** Throughout the proofs we repeatedly use the following inequalities on the quantiles of noncentral $\chi^2$ random variables. These inequalities are due to Birgé (2001). For all $u \in \,]0,1[$, $z \geq 0$, $d \geq 1$,

$$(18) \qquad q_{z,d}(u) \leq z + d + 2\sqrt{(2z+d)\log(1/u)} + 2\log(1/u),$$

$$(19) \qquad q_{z,d}(1-u) \geq z + d - 2\sqrt{(2z+d)\log(1/u)}.$$

In the sequel, $\Pi_m$ for $m \in \mathcal{M}_n$ denotes the orthogonal projector onto $S_m$.

5.1. *Proof of Lemma* 2.1. For simplicity, let us take $\sigma^2 = 1$.
If $\mathcal{D} = 0$, then $\hat{f} = \Pi_S Y = 0$, and hence

$$(20) \qquad \mathbb{P}_{f,1}[\phi(Y) = 0, \|f - \hat{f}\| \geq \rho] = \mathbb{P}_{f,1}[\|Y\|^2 \leq q_{0,n}(\alpha), \|f\| \geq \rho].$$

If $\|f\| < \rho$ this probability equals 0. Otherwise, $\|f\| \geq \rho$. Since $\|Y\|^2$ is distributed as a $\chi^2$ with noncentrality parameter $\|f\|^2$ and $n$ degrees of freedom, it follows from the definition of $\rho$ that the right-hand side of (20) is not larger than $\beta$.

Now let $\mathcal{D} \neq 0$. For all $f \in \mathbb{R}^n$, note that $\|\Pi_S \varepsilon\|^2$ and $\|Y - \Pi_S Y\|^2 = \|f - \Pi_S f + \varepsilon - \Pi_S \varepsilon\|^2$ are independent random variables. By setting $z = \|f - \Pi_S f\|^2$, we deduce

$$\mathbb{P}_{f,1}[\phi(Y) = 0, \|f - \hat{f}\| \geq \rho]$$

$$= \mathbb{P}_{f,1}[\|Y - \Pi_S Y\|^2 \leq q_{0,n-\mathcal{D}}(\alpha), \|f - \Pi_S f\|^2 + \|\Pi_S \varepsilon\|^2 \geq \rho^2]$$
$$= \chi^2_{z,n-\mathcal{D}}(q_{0,n-\mathcal{D}}(\alpha))(1 - \chi^2_{0,\mathcal{D}}(\rho^2 - z)).$$

If $\chi^2_{z,n-\mathcal{D}}(q_{0,n-\mathcal{D}}(\alpha)) \leq \beta$, then the result is established. Otherwise $z \in \mathcal{Z}$ and, by definition of $\rho$,

$$\rho^2 - z \geq q_{0,\mathcal{D}}\left(\frac{\beta}{\chi^2_{z,n-\mathcal{D}}(q_{0,n-\mathcal{D}}(\alpha))}\right),$$

which leads to

$$(1 - \chi^2_{0,\mathcal{D}}(\rho^2 - z)) \leq \frac{\beta}{\chi^2_{z,n-\mathcal{D}}(q_{0,n-\mathcal{D}}(\alpha))}$$

and the result follows.

5.2. *Proof of Theorems* 2.1 *and* 3.2. Theorem 2.1 being a straightforward consequence of Theorem 3.2 by taking $\eta = 0$, we only prove Theorem 3.2.

Let us first prove (17). The result is clear for $m = n$ as by definition $\hat{\rho} \leq \rho_n$. Let us fix some $m \in \mathcal{M}_n \setminus \{n\}$. We derive from the definition of $\hat{\rho}$ that

$$\mathbb{P}_{f,\sigma}[\hat{\rho} > \rho_m] \leq \mathbb{P}_{f,\sigma}[m \notin \mathcal{A}]$$
$$= \mathbb{P}_{f,\sigma}[\|Y - \hat{f}_m\|^2 > q_{0,N_m}(\alpha)\tau^2]$$
$$\leq \mathbb{P}_{f,\sigma}[\|Y - \hat{f}_m\|^2 > q_{0,N_m}(\alpha)\sigma^2],$$

as $\tau \geq \sigma$. We conclude by noting that, for $f \in S_m$, $\|Y - \hat{f}_m\|^2/\sigma^2$ is distributed as a $\chi^2$ with $N_m$ degrees of freedom.

We shall now show something that is stronger than (10), namely that

$$\mathbb{P}_{f,\sigma}\left[f \notin \bigcap_{m \in \mathcal{A}} \mathcal{B}(\hat{f}_m, \rho_m)\right] \leq \beta.$$

For all $f \in \mathbb{R}^n$,

$$\mathbb{P}_{f,\sigma}\left[f \notin \bigcap_{m \in \mathcal{A}} \mathcal{B}(\hat{f}_m, \rho_m)\right]$$
$$= \mathbb{P}_{f,\sigma}[\exists m \in \mathcal{A}, \|f - \hat{f}_m\| > \rho_m]$$
$$\leq \sum_{m \in \mathcal{M}_n} \mathbb{P}_{f,\sigma}[\|f - \hat{f}_m\| > \rho_m, \hat{m} \in \mathcal{A}]$$
$$= \sum_{m \in \mathcal{M}_n} \mathbb{P}_{f,\sigma}[\|f - \hat{f}_m\| > \rho_m, \|Y - \hat{f}_m\|^2 \leq q_{0,N_m}(\alpha)\tau^2].$$

Since $\sum_{m \in \mathcal{M}_n} \beta_m = \beta$, it is enough to prove that, for each $m \in \mathcal{M}_n$, the probability

$$\mathbb{P}_{f,\sigma}(m) = \mathbb{P}_{f,\sigma}[\|f - \hat{f}_m\| > \rho_m, \|Y - \hat{f}_m\|^2 \leq q_{0,N_m}(\alpha)\tau^2]$$

is not greater than $\beta_m$.

If $m = n$, this is clear since $Y = \hat{f}_n$ and, for $\tau^2 \geq \sigma^2$,

$$\mathbb{P}_{f,\sigma}(n) = \mathbb{P}_{f,\sigma}[\sigma^2 \|\varepsilon\|^2 > q_{0,n}(\beta_n)\tau^2] \leq \beta_n.$$

Let us now prove the inequality when $\mathcal{D}_m = 0$. In this case $\hat{f}_m = 0$. If $\|f\| \leq \rho_m$, we have $\mathbb{P}_{f,\sigma}(m) = 0$ and thus the inequality is true. Otherwise $\|f\| > \rho_m$ and as, for all $u > 0$ $z \to \chi^2_{z,n}(u)$ is nondecreasing with $z$ we get, by definition of $\rho_m$,

$$\mathbb{P}_{f,\sigma}(m) = \chi^2_{\|f\|^2/\sigma^2, n}(q_{0,n}(\alpha)\tau^2/\sigma^2)$$

$$\leq \chi^2_{\rho_m^2/\sigma^2, n}(q_{0,n}(\alpha)\tau^2/\sigma^2) \leq \beta_m.$$

Let us now fix some $m \in \mathcal{M}_n \setminus \{n\}$ such that $\mathcal{D}_m \neq 0$ and set $z = \|f - \Pi_m f\|^2/\sigma^2$. Note that the random variables

$$\frac{\|f - \hat{f}_m\|^2}{\sigma^2} = \frac{\|f - \Pi_m f + \sigma\Pi_m\varepsilon\|^2}{\sigma^2} = z + \|\Pi_m\varepsilon\|^2$$

and

$$\frac{\|Y - \hat{f}_m\|^2}{\sigma^2} = \frac{\|f - \Pi_m f + \sigma(\varepsilon - \Pi_m\varepsilon)\|^2}{\sigma^2}$$

are independent and that the second one is distributed as a noncentral $\chi^2$ with noncentrality parameter $z$ and $N_m$ degrees of freedom. Therefore, we get

$$(21) \qquad \mathbb{P}_{f,\sigma}(m) = \left(1 - \chi^2_{0,\mathcal{D}_m}\left(\frac{\rho_m^2}{\sigma^2} - z\right)\right)\chi^2_{z,N_m}\left(q_{0,N_m}(\alpha)\frac{\tau^2}{\sigma^2}\right).$$

We deduce from the definition of $\rho_m$ that, for all $\sigma \in I$ and $z \geq 0$, the right-hand side of (21) is not larger than $\beta_m$, which leads to the result.

5.3. *Proof of Theorem* 3.1. The principle of the proof leading to the lower bounds on the $r_m$'s is due to Lepski. However, the following nonasymptotic inequalities are to our knowledge new. In the sequel we set $N_m = n - \mathcal{D}_m$. Let us now fix some $m \in \mathcal{M}_n$; we divide the proof into consecutive claims.

CLAIM 1. *If $\alpha + \beta < 1 - \exp(-1/36)$, then*

$$r_m^2 \geq \left(\frac{\mathcal{D}_m}{27} - \sqrt{\mathcal{L}_1 \mathcal{D}_m}\right)\tau^2,$$

*where $\mathcal{L}_1 = -4\log(1 - \alpha - \beta)/81$.*

Note that the claim is clear when $\mathcal{D}_m = 0$; we shall thus restrict ourselves to the case $\mathcal{D}_m \geq 1$. The proof relies on two lemmas. In the first one, we show that, under the assumption of Theorem 3.1, with probability close to 1 the Euclidean distance between $f \in S_m$ and its estimator $\tilde{f}$ is not greater than $r_m$.

LEMMA 5.1. *Let the pair $(\tilde{f}, \tilde{r})$ satisfy the assumption of Theorem 3.1. Then, for all $m \in \mathcal{M}_n$, $f \in S_m$ and $\sigma \in I$,*

(22) $$\mathbb{P}_{f,\sigma}[\|f - \tilde{f}\| > r_m] \leq \alpha + \beta.$$

PROOF. For all $f \in S_m$,

$$\mathbb{P}_{f,\sigma}[\|f - \tilde{f}\| > r_m]$$
$$\leq \mathbb{P}_{f,\sigma}[\|f - \tilde{f}\| > r_m, r_m \geq \tilde{r}] + \mathbb{P}_{f,\sigma}[\|f - \tilde{f}\| > r_m, \tilde{r} > r_m]$$
$$\leq \mathbb{P}_{f,\sigma}[\|f - \tilde{f}\| > \tilde{r}] + \mathbb{P}_{f,\sigma}[\tilde{r} > r_m]$$

and we conclude thanks to (13) and (14). □

The second lemma shows that such a property of the estimator $\tilde{f}$ is possible only if $r_m$ is large enough.

LEMMA 5.2. *Let $S$ be a linear subspace of $\mathbb{R}^n$ of dimension $\mathcal{D} \geq 1$ and $\delta$ a positive number such that $\delta < 1 - \exp[-\mathcal{D}/36]$. If $\tilde{f}$ is an estimator of $f$ in (1) which satisfies, for all $f \in S$,*

(23) $$\mathbb{P}_{f,\sigma}[\|f - \tilde{f}\| > v_{\mathcal{D}}(\delta)] \leq \delta,$$

*then*

$$v_{\mathcal{D}}^2(\delta) \geq \left(\frac{\mathcal{D}}{27} - \frac{2}{9}\sqrt{\mathcal{D}\log(1/(1-\delta))}\right)\sigma^2.$$

In light of Lemma 5.1, the claim derives from Lemma 5.2 by taking $S = S_m$, $\delta = \alpha + \beta$ and $\sigma = \tau$. Let us now turn to the proof of Lemma 5.2.

PROOF OF LEMMA 5.2. The Gaussian law being invariant by orthogonal transformation, with no loss of generality, we assume that $S$ is the linear span generated by $e_1, \ldots, e_{\mathcal{D}}$, the $\mathcal{D}$ first vectors of the canonical basis of $\mathbb{R}^n$. Moreover, by homogeneity, we assume that $\sigma^2 = 1$. Let $v(\delta)$ be some positive number satisfying

(24) $$v^2(\delta) < \frac{\mathcal{D}}{27} - \frac{2}{9}\sqrt{-\mathcal{D}\log(1-\delta)}.$$

Note that the right-hand side of (24) is positive for $\delta \leq 1 - \exp[-\mathcal{D}/36]$. We prove Lemma 5.2 by showing that, for all estimators $\tilde{f}$ with values in $\mathbb{R}^n$,
$$\inf_{f \in S} \mathbb{P}_{f,1}[\|f - \tilde{f}\|^2 \leq v^2(\delta)] < 1 - \delta.$$

Let $\xi_1, \ldots, \xi_\mathcal{D}$ be Rademacher random variables (i.e., $\mathbb{P}[\xi_i = \pm 1] = 1/2$) which are independent of $Y$ and set $f(\xi) = \lambda \sum_{i=1}^{\mathcal{D}} \xi_i e_i$, where $\lambda$ denotes some positive number to be chosen later on. Using that
$$\frac{d\mathbb{P}_{f(\xi),1}}{d\mathbb{P}_{0,1}}(y) = \exp\left(-\frac{\lambda^2 \mathcal{D}}{2} + \lambda \sum_{i=1}^{\mathcal{D}} \xi_i y_i\right)$$
and the fact that $f(\xi) \in S$, we have
$$\inf_{f \in S} \mathbb{P}_{f,1}[\|f - \tilde{f}\|^2 \leq v^2(\delta)]$$
$$\leq \mathbb{P}_{f(\xi),1}\left[\sum_{i=1}^{\mathcal{D}}(\lambda \xi_i - \tilde{f}_i)^2 \leq v^2(\delta)\right]$$
$$= \mathbb{E}_{0,1}\left[\mathbf{1}\left\{\sum_{i=1}^{\mathcal{D}}(\lambda \xi_i - \tilde{f}_i(Y))^2 \leq v^2(\delta)\right\}\right.$$
$$\left. \times \exp\left(-\lambda^2 \mathcal{D}/2 + \lambda \sum_{i=1}^{\mathcal{D}} \xi_i Y_i\right)\right].$$

Note that $\tilde{f} = \tilde{f}(Y)$ satisfies
$$\sum_{i=1}^{\mathcal{D}}(\lambda \xi_i - \tilde{f}_i)^2 \geq \lambda^2 \sum_{i=1}^{\mathcal{D}} \mathbf{1}\{\xi_i \tilde{f}_i(Y) \leq 0\}$$
and thus, setting
$$N(\xi, \tilde{f}) = \lambda^2 \sum_{i=1}^{\mathcal{D}} \mathbf{1}\{\xi_i \tilde{f}_i(Y) \leq 0\},$$
we derive
$$\inf_{f \in S} \mathbb{P}_{f,\sigma}[\|f - \tilde{f}\|^2 \leq v^2(\delta)]$$
$$\leq \mathbb{E}_{0,1}\left[\mathbf{1}\{N(\xi, \tilde{f}) \leq v^2(\delta)\} \exp\left(-\lambda^2 \mathcal{D}/2 + \lambda \sum_{i=1}^{\mathcal{D}} \xi_i Y_i\right)\right].$$

By averaging with respect to $\xi$ and using Fubini's theorem we get
$$\inf_{f \in S} \mathbb{P}_{f,\sigma}[\|f - \tilde{f}\|^2 \leq v^2(\delta)]$$
(25)
$$\leq e^{-\lambda^2 \mathcal{D}/2} \mathbb{E}_{0,1}\left[\mathbb{E}_\xi\left[\mathbf{1}\{N(\xi, \tilde{f}) \leq v^2(\delta)\} \exp\left(\lambda \sum_{i=1}^{\mathcal{D}} \xi_i Y_i\right)\right]\right].$$

By the Cauchy–Schwarz inequality we have

$$\mathbb{E}_\xi^2\left[\mathbf{1}\{N(\xi,\tilde{f}) \leq v^2(\delta)\} \exp\left(\lambda \sum_{i=1}^{\mathcal{D}} \xi_i Y_i\right)\right]$$

$$\leq \mathbb{P}_\xi[N(\xi,\tilde{f}) \leq v^2(\delta)] \mathbb{E}_\xi\left[\exp\left(2\lambda \sum_{i=1}^{\mathcal{D}} \xi_i Y_i\right)\right]$$

$$= \mathbb{P}_\xi[N(\xi,\tilde{f}) \leq v^2(\delta)] \prod_{i=1}^{\mathcal{D}} \cosh(2\lambda Y_i),$$

which together with (25) gives

(26)
$$\inf_{f\in S} \mathbb{P}_{f,\sigma}[\|f-\tilde{f}\|^2 \leq v^2(\delta)]$$
$$\leq e^{-\lambda^2 \mathcal{D}/2} \mathbb{E}_{0,1}\left[\mathbb{P}_\xi^{1/2}[N(\xi,\tilde{f}) \leq v^2(\delta)] \prod_{i=1}^{\mathcal{D}} \cosh^{1/2}(2\lambda Y_i)\right].$$

Conditionally on $Y$, the random variable $N(\xi,\tilde{f})/\lambda^2$ is a sum of $\mathcal{D}$ independent random variables with values in $\{0,1\}$. Thus by Hoeffding's inequality we obtain that, for all $t \geq 0$,

$$\mathbb{P}_\xi[N(\xi,\tilde{f}) \leq \mathbb{E}_\xi[N(\xi,\tilde{f})] - \lambda^2\sqrt{\mathcal{D}t}] \leq e^{-2t}.$$

Taking $t = \lambda^2\mathcal{D}/2 - \log(1-\delta)$ and noting that $\mathbb{E}_\xi[N(\xi,\tilde{f})] \geq \lambda^2 \mathcal{D}/2$ we get from (24) that

$$\mathbb{E}_\xi[N(\xi,\tilde{f})] - \lambda^2\sqrt{\mathcal{D}t} \geq \lambda^2\left(\frac{\mathcal{D}}{2} - \sqrt{\frac{\lambda^2\mathcal{D}^2}{2} - \mathcal{D}\log(1-\delta)}\right)$$

$$\geq \left(\frac{\lambda^2}{2} - \frac{\lambda^3}{\sqrt{2}}\right)\mathcal{D} - \lambda^2\sqrt{-\mathcal{D}\log(1-\delta)}$$

and thus, for $\lambda = \sqrt{2}/3$,

$$\mathbb{E}[N(\xi,\tilde{f})] - \lambda^2\sqrt{\mathcal{D}t} \geq v^2(\delta).$$

Consequently,

$$\mathbb{P}_\xi^{1/2}[N(\xi,\tilde{f}) \leq v^2(\delta)] \leq e^{-t} = (1-\delta)e^{-\lambda^2 \mathcal{D}/2}.$$

Now using that

$$\mathbb{E}_{0,1}\left[\prod_{i=1}^{\mathcal{D}} \cosh^{1/2}(2\lambda Y_i)\right] = \prod_{i=1}^{\mathcal{D}} \mathbb{E}_{0,1}[\cosh^{1/2}(2\lambda Y_i)]$$
$$< \mathbb{E}_{0,1}^{\mathcal{D}/2}[\cosh(2\lambda Y_1)]$$
$$= \exp[\lambda^2 \mathcal{D}],$$

we derive from (26) that
$$\inf_{f \in S} \mathbb{P}_{f,\sigma}[\|f - \tilde{f}\|^2 \leq v^2(\delta)] < 1 - \delta,$$

which concludes the proof. □

CLAIM 2. *If $\alpha + 2\beta \leq 1 - \exp(-1/4)$, then*
$$(27) \qquad 9r_m^2 \geq \max\{\sqrt{\mathcal{L}_2 N_m}, (N_m - 2\sqrt{\mathcal{L}_3 N_m})\eta\}\tau^2,$$
*with $\mathcal{L}_2 = 2\log(1 + 4(1 - \alpha - 2\beta)^2)$ and $\mathcal{L}_3 = -\log(1 - \alpha - 2\beta)$.*

The claim is clear when $N_m = 0$; thus we only consider the case where $N_m \geq 1$. Again, the proof relies on two lemmas. The first one shows that if the pair $(\tilde{f}, \tilde{r})$ satisfies the assumptions of Theorem 3.1, then it is possible to build a level $(\alpha + \beta)$-test of "$f \in S_m$" against "$f \in \mathbb{R}^n \setminus S_m$" which achieves the power $1 - \beta$ on the complement of a ball of radius $3r_m$. Namely, the following holds:

LEMMA 5.3. *Let $(\tilde{f}, \tilde{r})$ be a pair of random variables with values in $\mathbb{R}^n \times \mathbb{R}_+$ satisfying the assumptions of Theorem 3.1. The test of hypothesis "$f \in S_m$" against the alternative "$f \notin S_m$" associated with the critical region*
$$(28) \qquad \mathcal{R} = \{\tilde{r} > r_m\} \cup \{\|\tilde{f} - \Pi_m \tilde{f}\| > 2\tilde{r}\}$$
*has the following properties: for all $\sigma \in I$,*
$$(29) \qquad \sup_{f \in S_m} \mathbb{P}_{f,\sigma}[\mathcal{R}] \leq \alpha + \beta,$$
*and for all $f$ satisfying $\|f - \Pi_m f\| > 3r_m$,*
$$(30) \qquad \mathbb{P}_{f,\sigma}[\mathcal{R}] \geq 1 - \beta.$$

PROOF. Let us show (29). First note that, for all $f \in S_m$,
$$(31) \qquad \begin{aligned} \|\tilde{f} - \Pi_m \tilde{f}\| &\leq \|f - \tilde{f}\| + \|f - \Pi_m \tilde{f}\| \\ &\leq 2\|f - \tilde{f}\|. \end{aligned}$$
By (13), (14) and (31), for all $f \in S_m$ we have
$$\mathbb{P}_{f,\sigma}[\mathcal{R}] \leq \mathbb{P}_{f,\sigma}[\tilde{r} > r_m]$$
$$+ \mathbb{P}_{f,\sigma}[\|\tilde{f} - \Pi_m \tilde{f}\| > 2\tilde{r}]$$
$$\leq \alpha + \mathbb{P}_{f,\sigma}[2\|f - \tilde{f}\| > 2\tilde{r}] \leq \alpha + \beta.$$

Let us now show (30). Let $f \in \mathbb{R}^n$ be such that $\|f - \Pi_m f\| \geq 3r_m$. Since
$$\|\tilde{f} - \Pi_m \tilde{f}\| \geq \|f - \Pi_m \tilde{f}\| - \|f - \tilde{f}\| \geq 3r_m - \|f - \tilde{f}\|,$$

we derive that

$$\begin{aligned}
\mathbb{P}_{f,\sigma}[\mathcal{R}^c] &= \mathbb{P}_{f,\sigma}[\|\tilde{f} - \Pi_m \tilde{f}\| \leq 2\tilde{r}, \tilde{r} \leq r_m] \\
&\leq \mathbb{P}_{f,\sigma}[\|\tilde{f} - \Pi_m \tilde{f}\| \leq 2r_m, \tilde{r} \leq r_m] \\
&\leq \mathbb{P}_{f,\sigma}[\|f - \tilde{f}\| \geq r_m, r_m \geq \tilde{r}] \\
&\leq \mathbb{P}_{f,\sigma}[\|f - \tilde{f}\| \geq \tilde{r}] \leq \beta. \qquad \square
\end{aligned}$$

We obtain the claim by proving that a test having the properties described in the previous lemma exists only if $r_m$ is large enough. The inequality

$$9r_m^2 \geq \sqrt{\mathcal{L}_2 N_m} \tau^2$$

derives from Baraud [(2002), Proposition 1]. For the second inequality,

$$9r_m^2 \geq (N_m - 2\sqrt{\mathcal{L}_3 N_m})\eta \tau^2,$$

we use the following lemma.

LEMMA 5.4. *Let $S$ be a linear subspace of $\mathbb{R}^n$ with $\dim(S) = \mathcal{D}$ (we set $N = n - \mathcal{D}$) and $\delta$ and $\beta$ be numbers satisfying $0 < \beta + \delta < 1 - \exp(-N/4)$. Let $\phi(Y)$ be a test function with values in $\{0, 1\}$ satisfying, for all $\sigma \in I$,*

(32) $$\sup_{f \in S} \mathbb{P}_{f,\sigma}[\phi(Y) = 1] \leq \delta,$$

*and for all $f \in \mathbb{R}^n$ such that $\|f - \Pi_S f\|^2 \geq \Delta(N, \beta)$,*

(33) $$\mathbb{P}_{f,\sigma}[\phi(Y) = 1] \geq 1 - \beta.$$

*Then*

$$\Delta(N, \beta) \geq (N - 2\sqrt{-N \log(1 - \beta - \delta)})\eta \tau^2.$$

By applying this lemma with $\delta = \alpha + \beta$, $S = S_m$ and $\mathcal{D} = \mathcal{D}_m$ and the test described in Lemma 5.3 we obtain the claim.

PROOF OF LEMMA 5.4. Let $\mathcal{F}$ be the set defined by

$$\mathcal{F} = \{f \in \mathbb{R}^n, \|\Pi_{S^\perp} f\|^2 \geq \Delta\},$$

where $\Delta$ denotes some positive number. To obtain the desired result it is enough to show that, for

$$\Delta < (N - 2\sqrt{-N \log(1 - \beta - \delta)})\eta \tau^2,$$

we have

(34) $$\inf_{\sigma \in I} \inf_{f \in \mathcal{F}} \mathbb{P}_{f,\sigma}[\phi(Y) = 1] < 1 - \beta.$$

Since the quantity $\sigma_* = \sqrt{1-\eta}\tau$ belongs to $I$, we have that, for all vectors $Z \in \mathbb{R}^n$,

$$\inf_{\sigma \in I} \inf_{f \in \mathcal{F}} \mathbb{P}_{f,\sigma}[\phi(Y) = 1]$$
$$\leq \mathbb{P}_{Z,\sigma_*}[\phi(Y) = 1]\mathbf{1}\{\|\Pi_{S^\perp}Z\|^2 \geq \Delta\} + \mathbf{1}\{\|\Pi_{S^\perp}Z\|^2 \leq \Delta\}.$$

By taking $Z$ as a random variable independent of $Y$ distributed as $\sqrt{\eta}\tau\varepsilon$, we obtain by averaging with respect to $Z$ that

$$\inf_{\sigma \in I} \inf_{f \in \mathcal{F}} \mathbb{P}_{f,\sigma}[\phi(Y) = 1] \leq \mathbb{E}[\mathbb{P}_{Z,\sigma_*}[\phi(Y) = 1]] + \mathbb{P}[\|\Pi_{S^\perp}Z\|^2 \leq \Delta].$$

For the first term of the right-hand side of this inequality, note that $\mathbb{E}[\mathbb{P}_{Z,\sigma_*}] = \mathbb{P}_{0,\tau}$. As $0 \in S$ and $\tau \in I$, we have

$$\mathbb{E}[\mathbb{P}_{Z,\sigma_*}[\phi(Y) = 1]] \leq \delta.$$

For the second term, note that our upper bound on $\Delta$ ensures that

$$\Delta < q_{0,N}(1 - \beta - \delta)\eta\tau^2$$

by using the lower bound on the quantiles of $\chi^2$ random variables (19). As the random variable $\|\Pi_{S^\perp}Z\|^2/(\eta\tau^2)$ is distributed as a $\chi^2(N)$, we get

$$\mathbb{P}[\|\Pi_{S^\perp}Z\|^2 \leq \Delta] < 1 - \beta - \delta,$$

which concludes the proof. $\square$

*Conclusion.* By gathering the inequalities of the two claims we get that, for some constant $C$ depending on $\alpha$ and $\beta$ only,

$$r_m^2 \geq C \max\{N_m\eta, \mathcal{D}_m, \sqrt{N_m}\}\tau^2.$$

Let us now prove (16). Let us fix some $f \in \mathbb{R}^n$. When $f = 0$, the result is clear by taking $S_m = \{0\}$. Then we deduce the result for general $f$ by arguing as follows. Let us consider the random variables $\tilde{f}_* = \tilde{f}(Y + f) + f$ and $\tilde{r}_* = \tilde{r}(Y + f)$. For all $g \in \mathbb{R}^n$ and $\sigma \in I$, we have that

$$\mathbb{P}_{g,\sigma}[g \in \mathcal{B}(\tilde{f}_*, \tilde{r}_*)] = \mathbb{P}_{g+f,\sigma}[g + f \in \mathcal{B}(\tilde{f}, \tilde{r})] \geq 1 - \beta.$$

Consequently, the pair of random variables $(\tilde{f}_*, \tilde{r}_*)$ satisfies (13) and thus, by taking $r_*(\alpha, 0) = r(\alpha, f)$ we derive that

$$r(\alpha, f) = r_*(\alpha, 0) \geq C \max\{\eta n, \sqrt{n}\}\tau^2.$$

5.4. *Proof of Propositions* 2.1 *and* 3.1. The result of the former proposition being a consequence of the latter by taking $\eta = 0$, we only prove Proposition 3.1. In the sequel we set $L_m = \log(1/\beta_m)$ and $L_\alpha = \log(1/\alpha)$. We distinguish three cases.

CASE $m = n$. We derive, from (18),
$$\rho_n^2 \leq (n + 2\sqrt{nL_n} + 2L_n)\tau^2,$$
which leads to the result.

CASE $\mathcal{D}_m \neq 0$, $m \neq n$. Let us fix $\sigma \in I$. Since for $z$ satisfying
$$\chi^2_{z,N_m}(q_{0,N_m}(\alpha)\tau^2/\sigma^2) \leq \beta_m$$
we have
$$(35) \quad z + q_{0,\mathcal{D}_m}\left(\frac{\beta_m}{\chi^2_{z,N_m}(q_{0,N_m}(\alpha)\tau^2/\sigma^2)} \wedge 1\right) = -\infty,$$

we bound from above the left-hand side of (35) for those $z$ satisfying
$$(36) \quad \chi^2_{z,N_m}(q_{0,N_m}(\alpha)\tau^2/\sigma^2) > \beta_m.$$

It follows from (19) that if $z$ satisfies (36), then
$$q_{0,N_m}(\alpha)\frac{\tau^2}{\sigma^2} \geq z + N_m - 2\sqrt{(2z + N_m)L_m}$$

and as we have
$$2\sqrt{(2z + N_m)L_m} \leq 2\sqrt{2zL_m} + 2\sqrt{N_m L_m} \leq \frac{z}{2} + 2\sqrt{N_m L_m} + 4L_m$$

and
$$q_{0,N_m}(\alpha) \leq N_m + 2\sqrt{N_m L_\alpha} + 2L_\alpha$$

from (18), we deduce that $z$ satisfies
$$(37) \quad \begin{aligned} z\sigma^2 &\leq \left(2\left(q_{0,N_m}(\alpha)\frac{\tau^2}{\sigma^2} - N_m\right) + 4\sqrt{N_m L_m} + 8L_m\right)\sigma^2 \\ &\leq \left(2N_m\eta + 4\sqrt{N_m}(\sqrt{L_m} + \sqrt{L_\alpha}) + 8L_m + 4L_\alpha\right)\tau^2. \end{aligned}$$

Thanks to (18) and the facts that $\chi^2_{z,N_m}(q_{0,N_m}(\alpha)\tau^2/\sigma^2) \leq 1$ and $\mathcal{D}_m \leq N_m$, we deduce that, for those $z$,
$$\begin{aligned} z\sigma^2 &+ q_{0,\mathcal{D}_m}\left(\frac{\beta_m}{\chi^2_{z,N_m}(q_{0,N_m}(\alpha)\tau^2/\sigma^2)} \wedge 1\right)\sigma^2 \\ &\leq \left(2N_m\eta + \mathcal{D}_m + 2\sqrt{N_m}(3\sqrt{L_m} + 2\sqrt{L_\alpha}) + 2(5L_m + 2L_\alpha)\right)\tau^2, \end{aligned}$$

and, consequently, that

$$\rho_m^2 \leq \left(2N_m\eta + \mathcal{D}_m + 2\sqrt{N_m}(3\sqrt{L_m} + 2\sqrt{L_\alpha}) + 2(5L_m + 2L_\alpha)\right)\tau^2.$$

The result follows as $N_m \leq n$.

CASE $\mathcal{D}_m = 0$. Arguing as above we have that for $x$ satisfying

$$x \geq \left(2N_m\eta + 4\sqrt{N_m}(\sqrt{L_m} + \sqrt{L_\alpha}) + 8L_m + 4L_\alpha\right)\tau^2$$

we have that, for all $\sigma \in I$,

$$\chi^2_{x/\sigma^2,n}(q_{0,n}(\alpha)\tau^2/\sigma^2) \leq \beta_m$$

and therefore, by definition of $\rho_m$,

$$\rho_m^2 \leq \left(2n\eta + 4\sqrt{n}(\sqrt{L_m} + \sqrt{L_\alpha}) + 8L_m + 4L_\alpha\right)\tau^2,$$

which leads to the result.

**Acknowledgments.** The author thanks the five referees for their suggestions that led to an improvement of the paper. The author is also grateful to Lucien Birgé for helpful discussions.

ECOLE NORMALE SUPÉRIEURE
DÉPARTEMENT DE MATHÉMATIQUES
ET APPLICATIONS
CNRS UMR 8553
45 RUE D'ULM
75230 PARIS CEDEX 05
FRANCE
E-MAIL: yannick.baraud@ens.fr